\documentclass[12pt]{amsart}

\usepackage{amsmath}
\usepackage{amssymb}
\usepackage{graphicx}
\usepackage{fullpage}
\usepackage{color}

\newcommand{\Yleft}{\includegraphics[height= 0.8 \baselineskip]{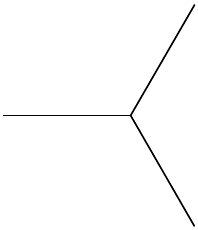}~ }
\newcommand{\Yright}{\includegraphics[height=0.8 \baselineskip]{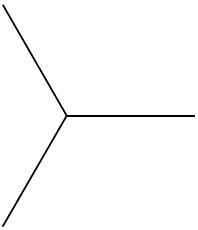}~ }

\newcommand{\bbZ}{\mathbb{Z}}
\newcommand{\bbC}{\mathbb{C}}
\newtheorem{theorem}{Theorem}
\newtheorem{lemma}[theorem]{Lemma}
\newtheorem{definition}[theorem]{Definition}
\newcommand{\pf}{\noindent \textbf{Proof. }}

\begin{document}

\title{Squishing dimers on the hexagon lattice}
\email{byoung@math.mcgill.ca}
\author{Ben Young}
\maketitle

\begin{abstract}
We describe an operation on dimer configurations on the hexagon lattice, called ``squishing'', and use this operation to explain some of the properties of dimer generating functions.
\end{abstract}

\section{Introduction}

In this paper, we will describe and use a novel technique called ``squishing'', which one applies to the dimer model on the regular honeycomb lattice.  We developed this technique as an attempt to verify a conjectured generating function which arises in algebraic geometry (specifically, in the Donaldson-Thomas theory of the orbifold $\bbC^3 / \bbZ_2 \times \bbZ_2$~\cite{YOUNG-BRYAN}).  Our attempt was only partially successful, and we were later able to compute the generating function by other means.  However, the technique is interesting in of itself, being one of relatively few dimer model techniques which exploits the self-similarity of the lattice at different scales.

We begin by describing the original motivation for this work.

\begin{definition}
A \emph{3D Young diagram} (or \emph{3D diagram}) $\pi$ is a subset of $(\bbZ_{\geq 0})^3$ such that if $(x,y,z) \in \pi$, then $(x',y',z') \in \pi$ whenever $x' \leq x$, $y' \leq y$, and $z' \leq z$.  
\end{definition}
3D Young diagrams are called \emph{boxed plane partitions} or \emph{3D partitions} elsewhere in the literature.  We refer to the points in $\pi$ as \emph{boxes} -- the point $(i,j,k)$ corresponds to the unit cube with vertices $\{(i \pm \frac{1}{2}, j \pm \frac{1}{2}, k \pm \frac{1}{2})  \}$.

We will be discussing the following generating functions:

\begin{definition} A \emph{weighting} of $(\bbZ_{\geq 0})^3$ is a map
\[
w:(\bbZ_{\geq 0})^3 \rightarrow \{p,q,r,s\},
\] 
where $p,q,r,s$ are formal indeterminates.   $w(p)$ is called \emph{weight} of the lattice point $p$;  the \emph{weight} of a 3D diagram is the product of the weights of the lattice points at the centers of all of its boxes.  The \emph{$w$-partition function} is then defined to be the formal sum
\[
Z_w = \sum_{\pi \text{ 3D diagram}} w(\pi).
\]
\end{definition}
The weightings we will be concerned with are the \emph{monochromatic} weighting $w_{\{1\}}$,
\[
(i,j,k) \mapsto p
\]
and the $\bbZ_2 \times \bbZ_2$ weighting $w_{\bbZ_2 \times \bbZ_2}$,
\[
(i,j,k) \mapsto \begin{cases} 
p & \text{if } i-k \equiv 0, j-k \equiv 0 \pmod{2}, \\
q & \text{if } i-k \equiv 1, j-k \equiv 0 \pmod{2},\\
r & \text{if } i-k \equiv 0, j-k \equiv 1 \pmod{2},\\
s & \text{if } i-k \equiv 1, j-k \equiv 1 \pmod{2},
\end{cases}
\]
along with various specializations of these weightings; we shall denote their partition functions $Z_{\{1\}}(p)$ and $Z_{\bbZ_2 \times \bbZ_2}(p,q,r,s)$, respectively.

It is a classical result~\cite{MACMAHON} that
\begin{equation}
\label{eqn:monochrome_function}
\sum_{\pi} w_{\{1\}}(\pi) = M(1,p)
\end{equation}
where we define
\[
M(a,z) = \prod_{n=1}^{\infty} \left(\frac{1}{1-az^n}\right)^n.
\]

Equation~\eqref{eqn:monochrome_function} also arises in algebraic geometry, essentially because 3D Young diagrams are the same as monomial ideals $I \subset \bbC^3[x,y,z]$ (simply read off the exponents of the elements of the coordinate ring $\bbC^3 / I$; these are the boxes of $\pi$).  As a result, ~\eqref{eqn:monochrome_function} is, in a certain sense, an invariant of $\bbC^3$; specifically, it is the \emph{Donaldson-Thomas partition function} for the space $\bbC^3$, up to a sign on $p$.
\cite{ORV}.

It is possible~\cite{YOUNG-BRYAN} to develop Donaldson-Thomas theory for orbifolds of $\bbC^3$ under the actions of certain finite abelian groups.  It turns out that for the group $\bbZ_2 \times \bbZ_2$, the partition function is given by

\begin{equation}
\label{eqn:z2z2_function}
Z_{\bbZ_2 \times \bbZ_2} =
\frac{M(1,Q)^4 \widetilde{M}(qr, Q) \widetilde{M}(qs, Q) \widetilde{M}(rs, Q)}
{\widetilde{M}(-q, Q) \widetilde{M}(-r, Q)\widetilde{M}(-s, Q)\widetilde{M}(-qrs, Q)   }
\end{equation}
where $Q=pqrs$ and $\widetilde{M}(a, z) = M(a,z) M(a^{-1},z)$. 

The curious identity \eqref{eqn:z2z2_function} was proven in~\cite{YOUNG-BRYAN} using vertex operators; it was conjectured earlier by Jim Bryan (based on the behaviour of related Donaldson-Thomas partition functions) and independently by Rick Kenyon (in an equivalent form, based on empirical computer work).

It is possible to check some of the properties of~\ref{eqn:z2z2_function} in an elementary manner.  For example, it should be the case that specializing $p=q=r=s$ should give
\[Z_{\bbZ_2 \times \bbZ_2}(p,p,p,p) = Z_{\{1\}}, \] and indeed~\eqref{eqn:z2z2_function} does satisfy this relation.  Another striking relation suggested by \eqref{eqn:z2z2_function} is:
\begin{equation}
\label{eqn:-1_specialization}
Z_{\bbZ_2 \times \bbZ_2}(p, -1, -1, -1) = M(1, Q)^2;
\end{equation}
however, the combinatorial reason for this is far less clear.  
This paper demonstrates~\eqref{eqn:-1_specialization} in an elementary manner, without relying on the theorems of~\cite{YOUNG-BRYAN}.  

\section{Matchings, and Weightings}

Our attack on ~\eqref{eqn:-1_specialization} immediately requires us to to encode the ``surface'' of a 3D Young diagram with dimers on the hexagon lattice.  Let us fix some terminology.

\begin{definition} Let $G = (V,E)$ be a graph.  A \emph{matching} or \emph{1-factor} of $G$ is a subgraph $M=(V,E')$ such that the degree of every vertex $v \in V$ is 1.  Equivalently, a matching of $G$ is a partition of the vertices of $G$ into a disjoint union of \emph{dimers}, or pairs of vertices joined by a single edge of $G$.
\end{definition}

\begin{figure}
\caption{A matching of $H_{5,4,3}$.}
\label{fig:matching}
\begin{center}
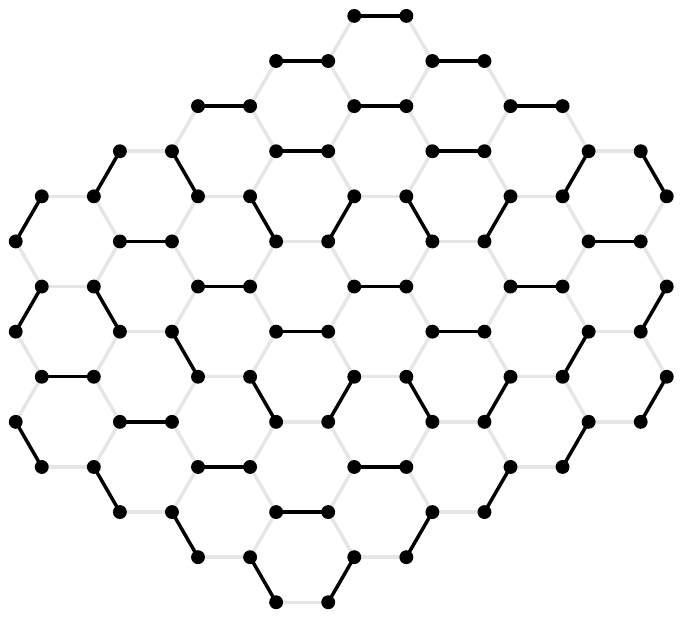
\end{center}
\end{figure}

For a general survey of results on the dimer model, see~\cite{dimer-survey}.

\begin{figure}
\caption{A 3D Young diagram viewed as a matching on $H_{3,3,3}$.}
\label{fig:matching_as_part}
\begin{center}
\makebox{
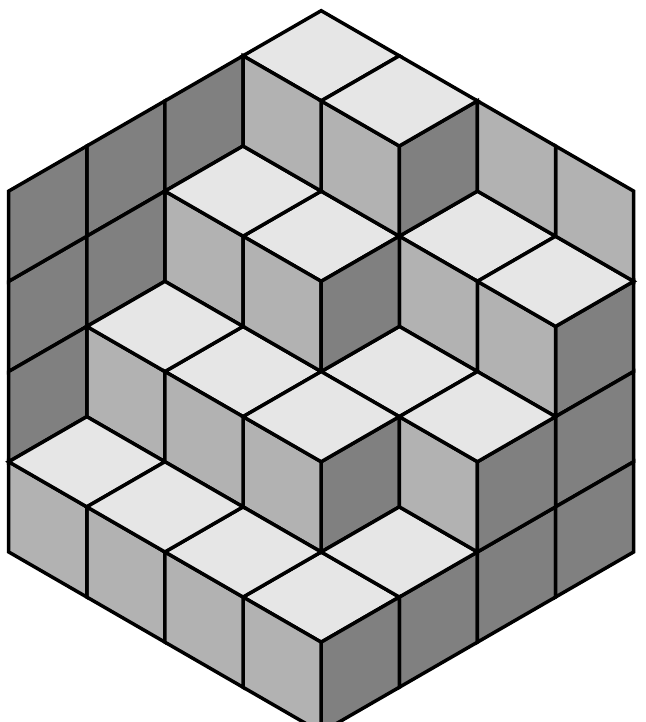
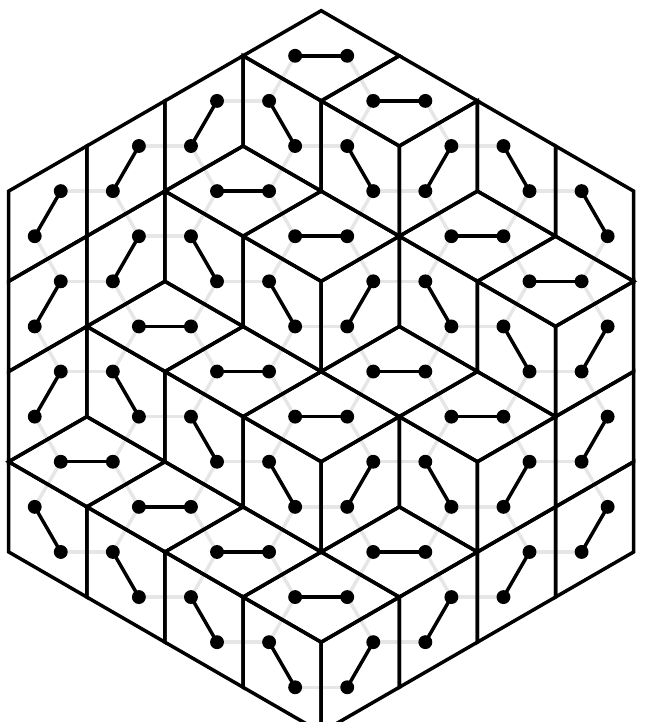
}
\end{center}
\end{figure}

We will be considering matchings on the semiregular hexagonal mesh of side
lengths $a,b,c,a,b,c$, denoted $H_{a,b,c}$ (see Figure~\ref{fig:matching} for a definition--by--picture).  Matchings on this graph are in
bijection with 3D Young diagram which are contained within an $a \times b
\times c$ box.  To see why, imagine viewing a 3D Young diagram from a distance(i.e., under isometric projection).  The faces of the 3D diagram are then rhombi, each of which is composed of two equilateral triangles (see Figure~\ref{fig:matching_as_part}).  The centers of these triangles fall at the vertices of $H_{a,b,c}$, so we get a matching by replacing each rhombus with the corresponding edge.  

\begin{figure}
\caption{A monochromatic weighting on $H_{4,4,4}$}
\label{fig:monochromatic_weighting}
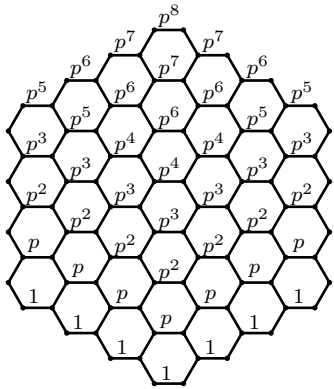
\end{figure}

A weighting of a graph assigns a monomial to each of the graph's edges.  Our first example, called a \emph{monochromatic weighting}, is the one depicted in figure~\ref{fig:monochromatic_weighting}, which assigns a weight of one to all non-horizontal edges, and weight $p^i$ to horizontal edges which have $i$ other horizontal edges directly below them.  We will adopt the convention that if an edge has no weight written beside it, then that edge has a weight of one.  

If our graph has a weighting, then the weight of a matching is the product of
all of the weights of the edges of the graph.  

\begin{figure}
\caption{The empty 3D Young diagram and its associated $w_p$-weighted matching}
\label{fig:emptypart}
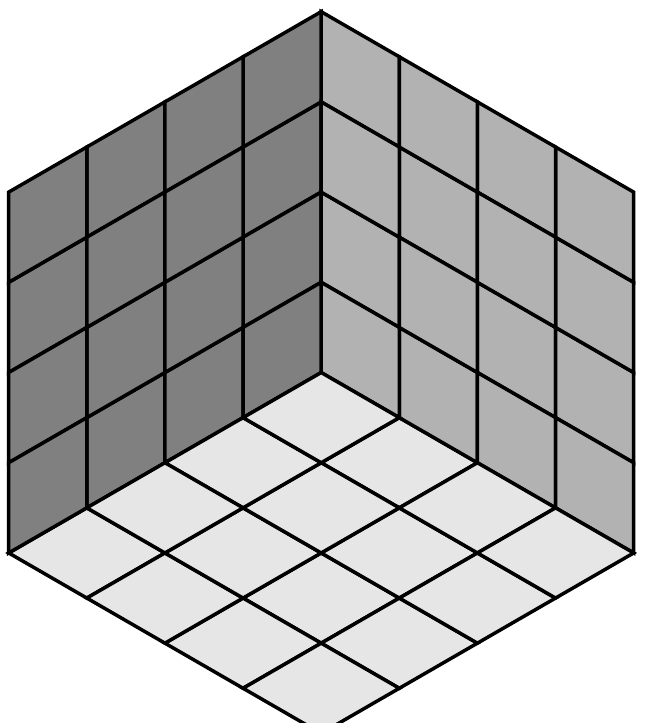
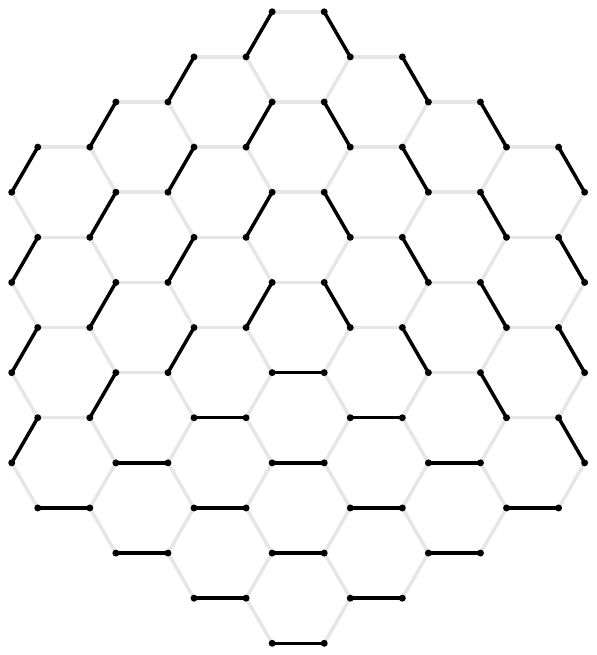
\end{figure}

Now, we have defined two monochromatic weightings -- one for 3D Young diagrams and one for graphs.  One might ask whether they are ``the same'': is the weight of a 3D diagram which fits inside an $i\times j\times k$ box the same as the weight of the associated matching of $H_{i,j,k}$?  Strictly speaking, the answer is no, because the weight of the empty 3D Young diagram should be zero, but the weight of the associated matching (see Figure~\ref{fig:emptypart}) is nonzero.

However, given a matching $M$, we can define its \emph{normalized weight} as the weight of $M$ divided by the weight of the empty 3D Young diagram.  Now it is easy to see that the normalized weight of $M$ is equal to the weight of the associated 3D diagram $\pi$. The proof is by induction on the number $n$ of boxes in $\pi$, the case $n=0$ being trivial.  One only needs to check that the operation of adding a box to $\pi$ has the effect of increasing both weights by $p$, which is easy to do.

\begin{figure}
\caption{A ``better'' monochromatic weighting, $w_p$,  where $p=t^3$.}
\label{fig:p_weighting}
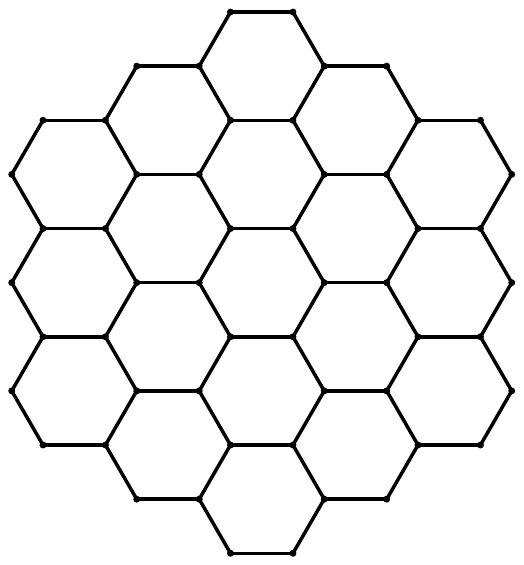
\end{figure}

There are many other weightings on the hexagon meshes whose normalized versions are equivalent to the monochromatic weight.  For example, one could rotate the weighted mesh by 120 degrees.   Indeed, for our purposes, the following weighting is superior (see Figure~\ref{fig:p_weighting}): we replace $p$ with $t$, and superimpose three copies of the old monochromatic weighting:  one normal, one rotated 120 degrees, and one weighted 240 degrees.  This weighting assigns each box the weight $t^3$, so substituting $p=t^3$ gives a new monochromatic weighting.

\begin{definition} The weighting described above is called $w_p$.  
\end{definition}
It is cumbersome to draw diagrams with $p^{1/3}$ as an edge weight, so throughout the paper, we shall use the convention that $p=t^3$when it is convenient.

\begin{figure}
\caption{The $qt$--, $rt$--, and $st$--weightings on $H_{4,4,4}$.}
\label{fig:pqrs_weightings}
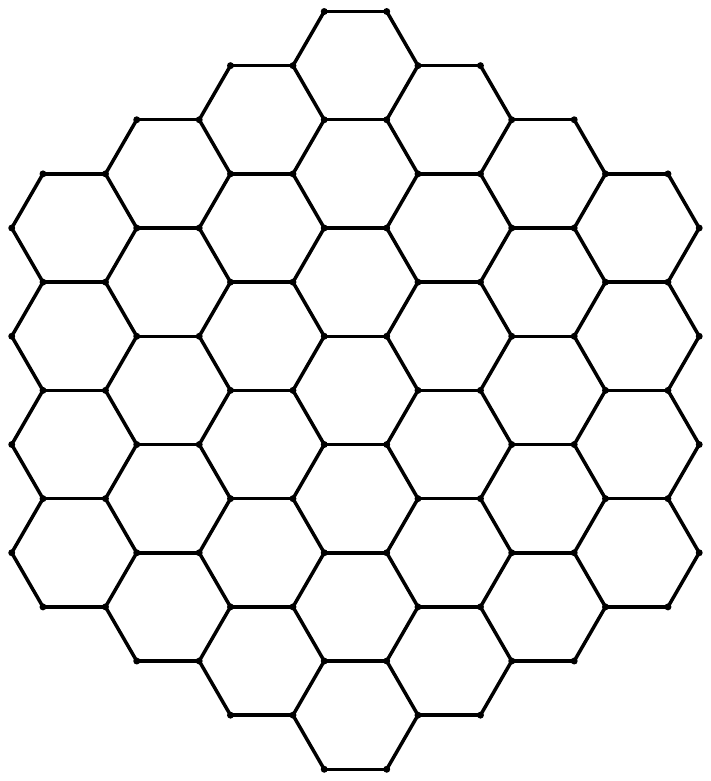
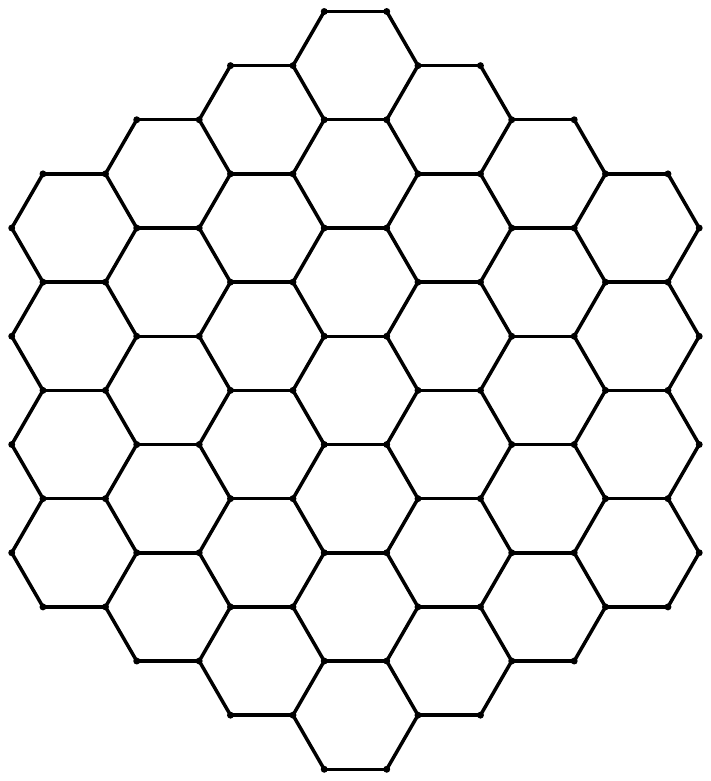

\vspace{-0.6in}

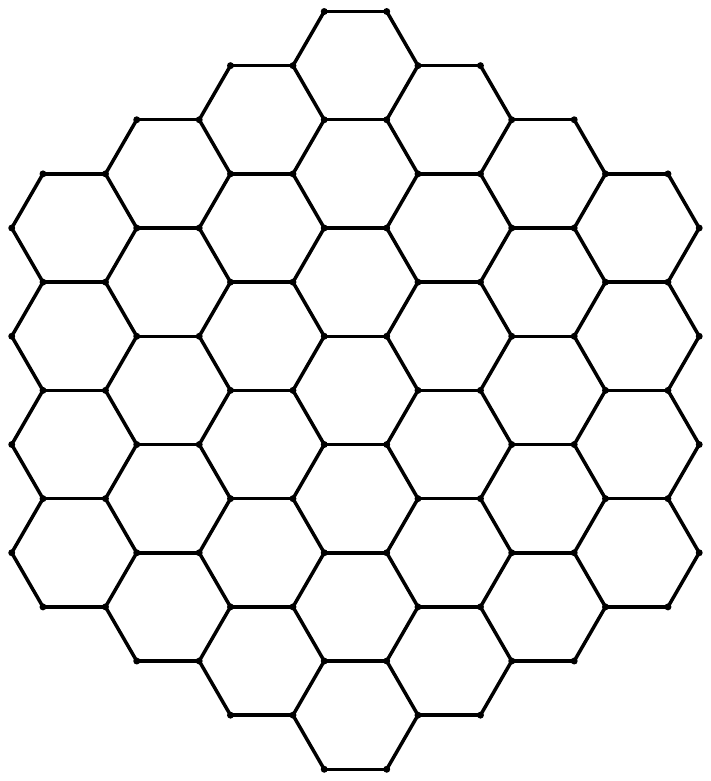
\end{figure}

We would next like to define a weighting of $H_{a,b,c}$ whose normalized version is equivalent to the $\bbZ_2 \times \bbZ_2$ weighting.  There are several ways of doing this, but for our purposes, the best way is to first define three weightings: the $qt$--, $rt$--, and $st$--weightings (see Figure~\ref{fig:pqrs_weightings}).   The $qt$--weighting is equivalent to the $\bbZ_2 \times \bbZ_2$ weighting under the specialization $p \mapsto t; r,s \mapsto 1$, and similarly for the other two weightings.

\begin{figure}
\caption{The $\bbZ_2 \times \bbZ_2$ weighting.}
\label{fig:z2z2_weighting}
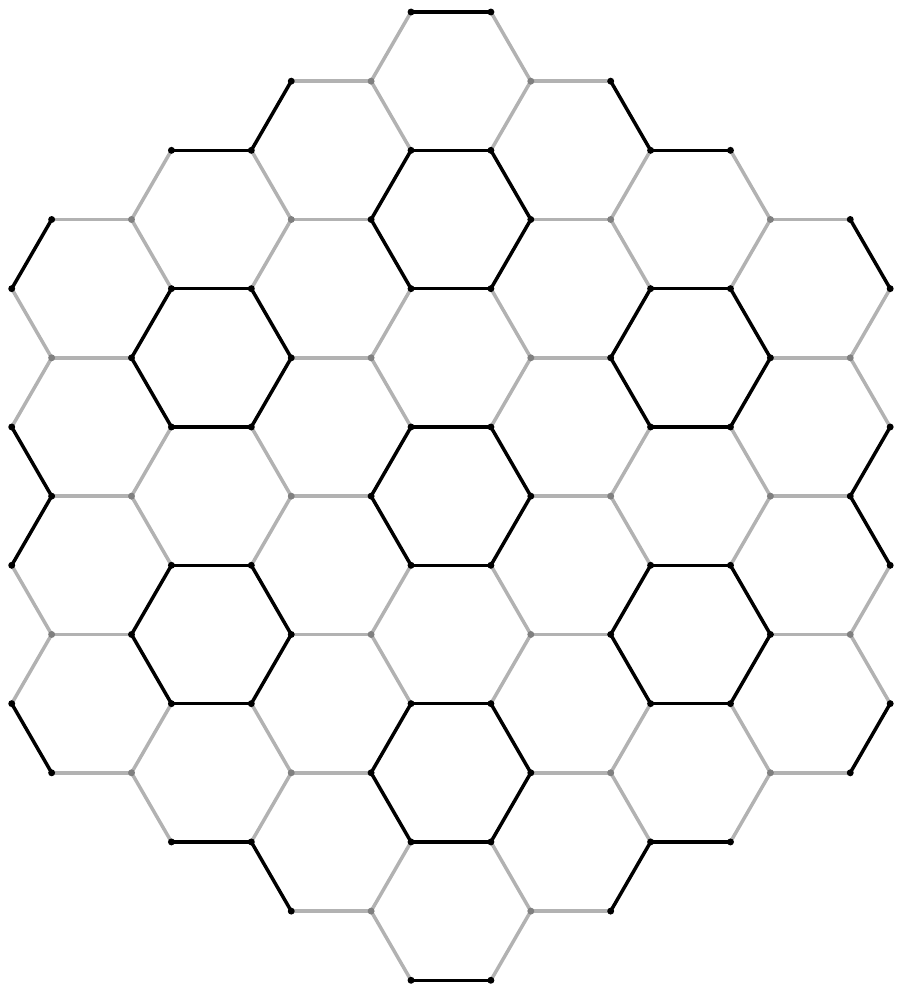
\end{figure}

Having done this, we construct the $\bbZ_2 \times \bbZ_2$ weighting by assigning each edge in $H_{a,b,c}$ the product of its $qt$,$rt$, and $st$--weights.  This weights boxes colored $q,r,s$ correctly, but each box in the $p$ position gets the weight $t^3$.  So specializing $t \mapsto p^{1/3}$ gives the $\bbZ_2 \times \bbZ_2$ weighting (see Figure~\ref{fig:z2z2_weighting}).  Observe that we have assigned weight 1 to all of the grey edges.

\begin{definition}
We call the weighting of Figure~\ref{fig:z2z2_weighting} $w_{p,q,r,s}$.
\end{definition}

\section{Overlaying pairs of matchings}

We were introduced to the ideas in this section by Kuo's beautiful paper on graphical condensation~\cite{Kuo}.  In the following, $G$ will always be a bipartite graph.

Suppose that we have two matchings $M_1, M_2$ of $G$.  If we overlay these two matchings on the vertex set $V$ of $G$, we have a multigraph $N$ in which each vertex has degree two, called a \emph{2-factor} of $G$.   This terminology is slightly nonstandard in graph theory -- elsewhere in the literature, a 2-factor is usually a collection of closed loops and isolated edges (not doubled).

If the edge $e$ occurs in both $M_1$ and $M_2$, then $e$ occurs as a doubled edge in $N$.  In this case, since the degree of both endpoints of $e$ is two, $e$ is a connected component in $N$.  Conversely, all doubled edges in $N$ must occur in both $M_1$ and $M_2$.  If we disregard the doubled edges, the rest of $N$ decomposes into a collection of disjoint closed paths.

Conversely, one can split a 2--factor into two one-factors:

\begin{lemma}
A 2-factor $N$ may be partitioned into an ordered pair of matchings $(M_1, M_2)$ in precisely $2^{\#\{\text{closed paths in }N\}}$ distinct ways.
\end{lemma}

\pf
Suppose we have a 2--factor $N$ of a bipartite graph $G$.  We may obtain two matchings $M_1$ and $M_2$ of $G$ as follows:  If $e$ is doubled in $N$, then place $e$ into both $M_1$ and $M_2$.  If $P$ is a closed path in $N$, select one of the edges in $P$ and place it into $M_1$.  Place the next edge in the path into $M_2$, and so forth.  Since $G$ is bipartite, the path $P$ is of even length, so each vertex in $P$ has degree 1 in both $M_1$ and $M_2$.

There are $2$ ways to divide $P$ between $M_1$ and $M_2$, so there are $2^{\#\{\text{closed paths in } N\}}$ pairs of matchings $M_1, M_2$ which correspond to $N$. $\hfill \square$

\begin{figure}
\caption{Overlaying two matchings on $H_{3,3,3}$.}
\label{fig: overlay}
\begin{center}
\makebox{
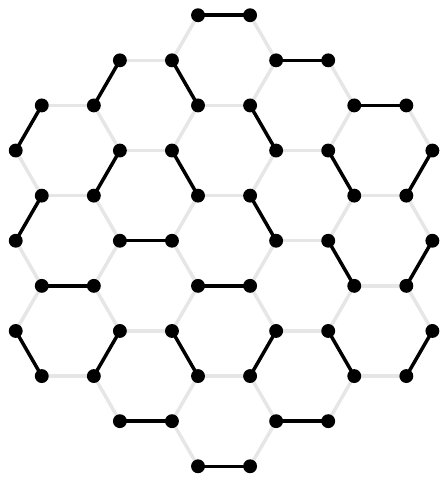
\negthickspace \negthickspace \negthickspace \negthickspace \negthickspace
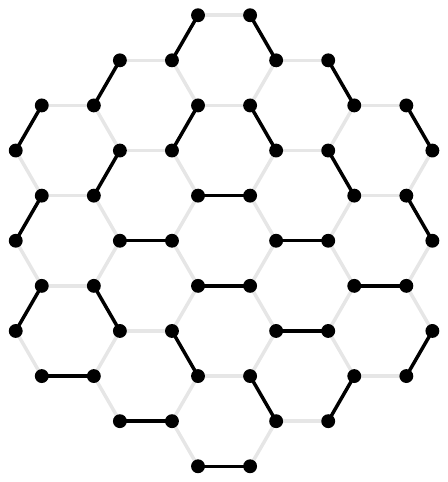
\negthickspace \negthickspace \negthickspace \negthickspace \negthickspace
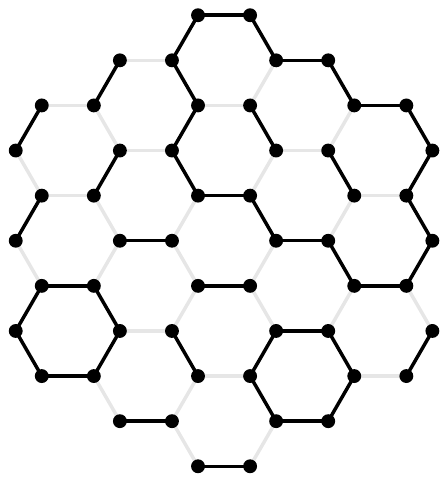
}
\end{center}
\end{figure}

\section{Squishing}

\begin{figure}
\caption{The hexagonal mesh $H_{4,4,4}$ being squished.}
\label{fig:squish}
\begin{center}
\makebox{
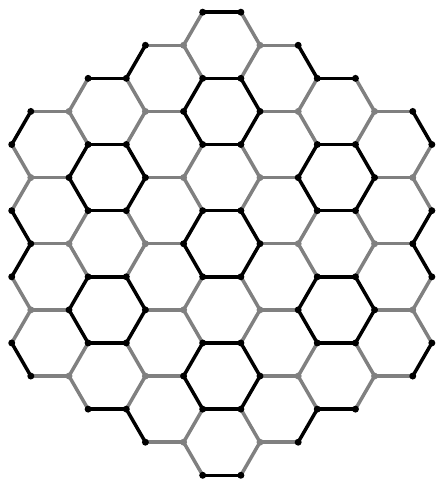
\negthickspace \negthickspace \negthickspace \negthickspace \negthickspace
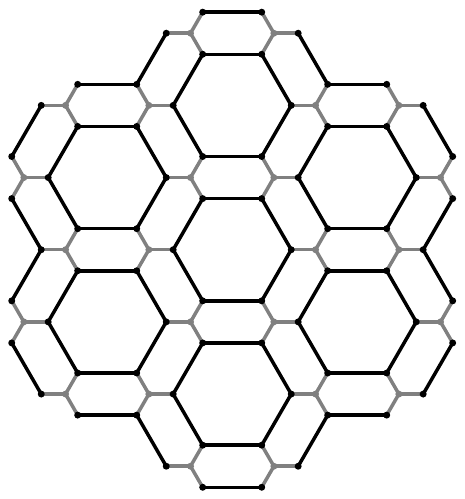
\negthickspace \negthickspace \negthickspace \negthickspace \negthickspace
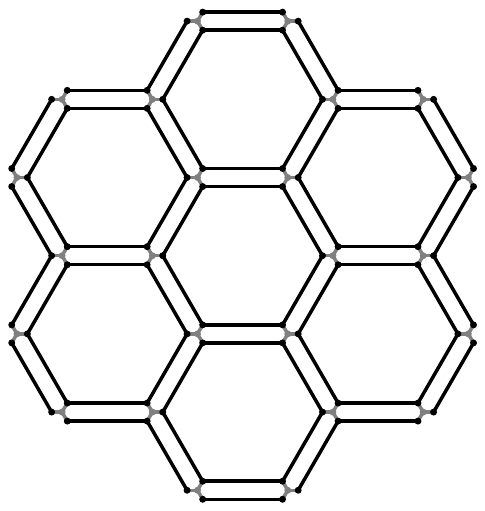
}
\end{center}
\end{figure}

Consider the hexagonal mesh with even side lengths, $H_{2a,2b,2c}$.  The
leftmost diagram in Figure~\ref{fig:squish} shows a picture of $H_{4,4,4}$, with
some of the edges colored grey.  The grey edges come in sets of three, all of
which are incident to one central vertex.  I shall call these three edges a
\emph{propeller}.

The other two diagrams in Figure~\ref{fig:squish} show what happens when the length of the edges in each propeller is decreased, while the other edges remain long.   We call this procedure \emph{squishing}.   Of course, the length that we choose to draw the edges in the graph has no bearing of the structure of the graph, but when the propellers are quite small, the graph of $H_{4,4,4}$ ``looks like'' the graph of $H_{2,2,2}$ with each edge doubled.  

We will denote this squishing operation by the symbol $\psi:$
\[
\psi:H_{2a,2b,2c} \setminus \{\text{propellers}\} \longrightarrow H_{a,b,c}
\]
where $\psi$ sends each edge to its image after squishing.  It is ofted useful to speak of squishing a set $E$ of edges of $H_{2a,2b,2c}$, and we shall also denote this operation by $\psi$:
\[
\psi(E) = \bigcup_{e \in E \setminus \{\text{propellers} \}} \psi(e).
\]

\begin{figure}
\caption{The ``unsquishing'' map $\varphi$}
\label{fig:maximal_unsquishing}
\begin{center}
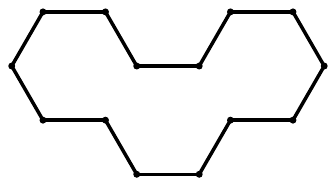
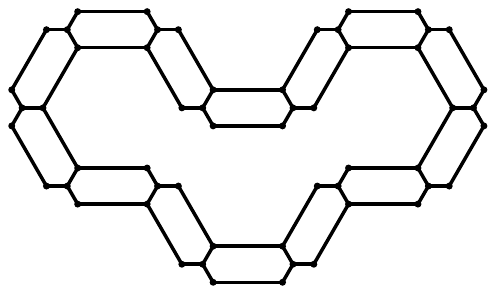
\end{center}
\end{figure}

Sometimes, given $E' \subset H_{a,b,c}$, we will need to look for sets of edges $E$ for which $\psi(E) = E'$.  Naturally, there are many such $E$, since $\psi$ ignores the propellers and is two-to-one on all other edges.  However, for a given $E' \subset H_{a,b,c}$, there is a ``most relevant'' preimage of $E'$ under $\psi$, which we shall call $\varphi(E') \subset H_{2a,2b,2c}$, defined as follows: 
\begin{align*}
E_0 & := \bigcup_{e \in E'} \phi^{-1}(e), \\
\varphi(E') &:= E_0 \cup \{ \text{all propellers adjacent to }E_0\}
\end{align*}
In other words, $\varphi(E')$ contains all edges which squish to edges of $E'$, as well as all propellers incident to those edges (See Figure~\ref{fig:maximal_unsquishing}).  It is clear that $\psi \circ \varphi(E') = E'$.

\section{Three lemmas on weightings of squished graphs}

\begin{figure}
\caption{A matching on $H_{4,4,4}$ being squished to a 2--factor on $H_{2,2,2}$.}
\begin{center}
\label{fig:squishmatching}
\makebox{
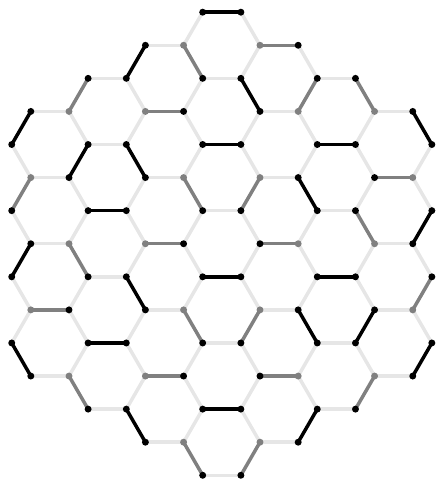
\negthickspace \negthickspace \negthickspace \negthickspace \negthickspace
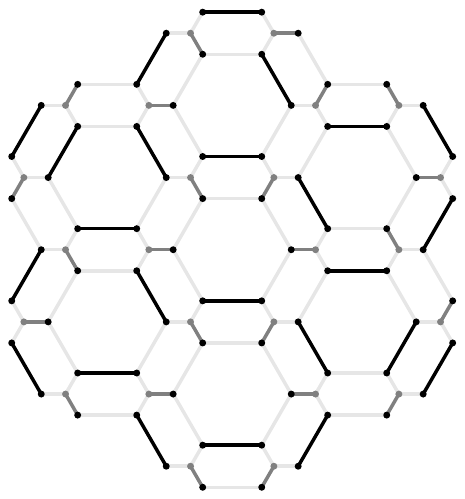
\negthickspace \negthickspace \negthickspace \negthickspace \negthickspace
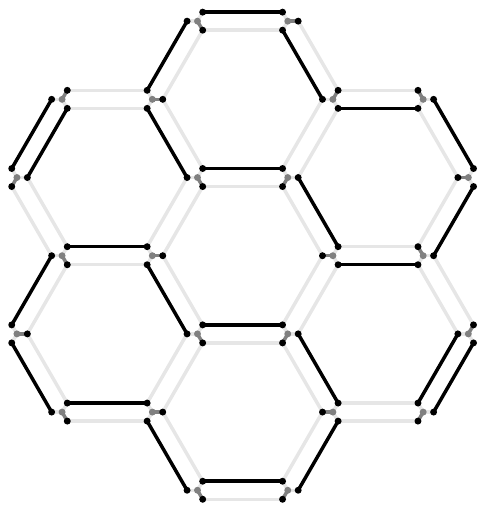
}
\end{center}
\end{figure}

When one draws a matching on the graph of $H_{4,4,4}$ before squishing, we get what looks very much like a 2-factor of $H_{2,2,2}$ (See Figure~\ref{fig:squishmatching}).  We will quantify precisely how this occurs, and use the effect to prove several facts about certain specializations of the $\bbZ_2 \times \bbZ_2$ partition function.

\begin{figure}
\caption{All possible configurations of edges around a propeller in a perfect matching.}
\begin{center}
\label{fig:propeller_possibilities}
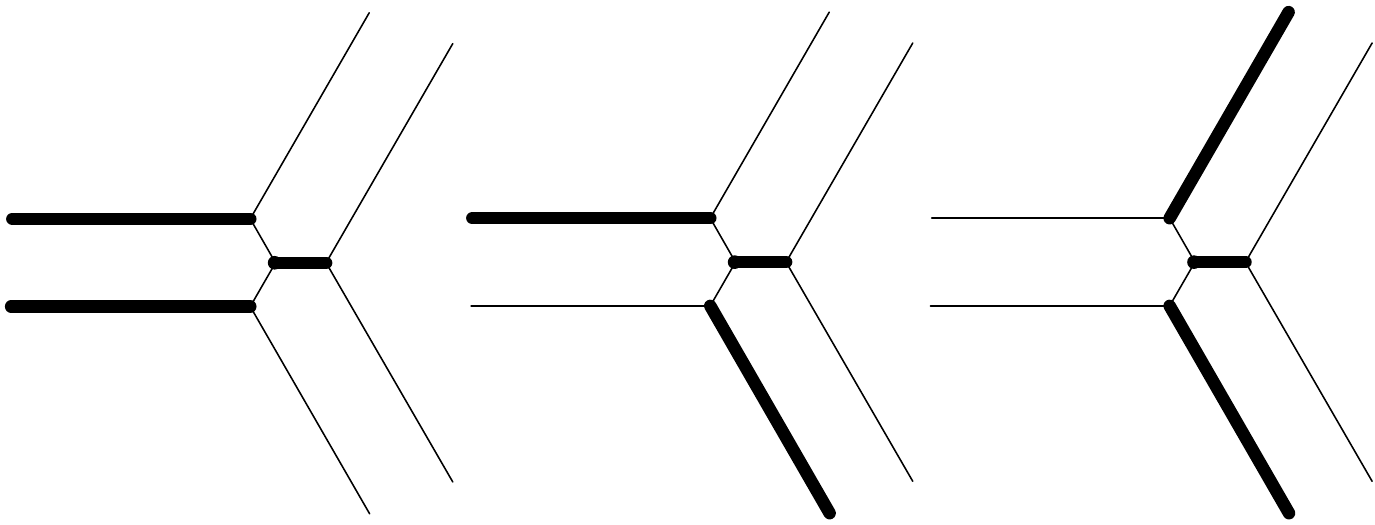
\end{center}
\end{figure}

Let us examine one of the propellers more closely, and determine what possible configurations of its edges can appear in a perfect matching $M$, up to symmetry.  (see Figure ~\ref{fig:propeller_possibilities}).   We label the central vertex $D$, and the other vertices $A,B,C$.  First, observe that preciesly one of the three "short" edges must be included in $M$, to give the central vertex degree one.  Suppose that it is the short horizontal edge $CD$.  In order for vertices $A$ and $B$ to have degree 1 as well, they must each have one incident long edge.  Up to symmetry, there are three ways for this to occur:
\begin{enumerate}
\item All three edges are horizontal;
\item The edge incident to $A$ is horizontal; the edge incident to $B$ is diagonal; or
\item The edges incident to both $A$ and $B$ are diagonal.
\end{enumerate}

If we now identify the vertices $A,B,C,D$ and ignore the short edges, we can see that case 1 gives rise to a doubled edge, and cases 2 and 3 give rise to a path through the vertex.   Therefore, collapsing all propellers in $H_{2a,2b,2c}$ to points transforms a matching on $H_{2a, 2b, 2c}$ into a 2-factor on  $H_{a,b,c}$.  In other words, the squishing operation $\psi$ induces a map, which we will call $\Psi$:
\[
\Psi: \{\text{1-factors of }H_{2a, 2b, 2c}\} \longrightarrow \{\text{2-factors of }H_{a,b,c}\}.
\]  
It is easy to see that $\Psi$ is surjective but not injective.  For example, $H_{2,2,2}$ has 20 perfect matchings, whereas $H_{1,1,1}$ has only 3 2-factors.  

It is often necessary to find the set of all matchings on $\varphi(E')$.  Strictly speaking, this set is $\left(\Psi |_{\varphi(E')} \right)^{-1}(E')$, but for simplicity of notation we shall just write this as $\Psi^{-1}(E')$. 

\begin{figure}
\caption{The $\bbZ_2 \times \bbZ_2$ weighting, under the specialization $t,q,r,s \mapsto -1$}
\label{fig:minus_1_weighting}
\begin{center}
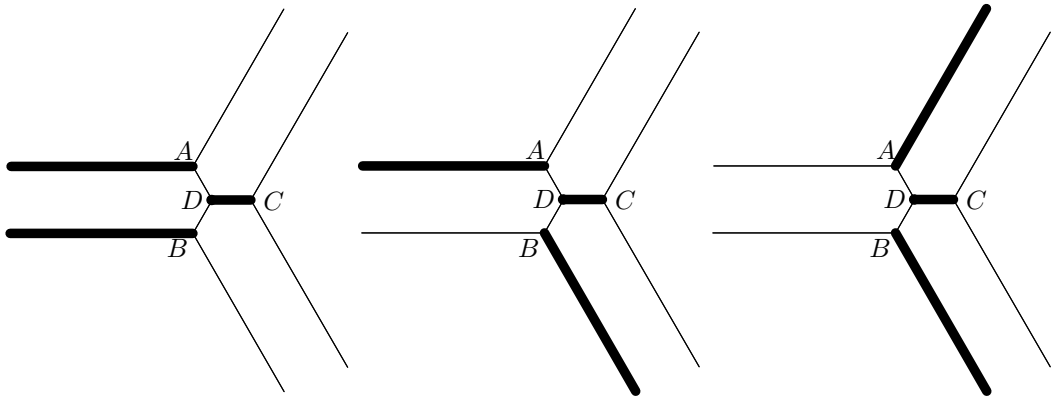
\end{center}
\end{figure}

We will need to use the $\bbZ_2 \times \bbZ_2$ weighting under the specialization $t,q,r,s \mapsto -1$, which we shall call $w_{-1}$ for short (see Figure~\ref{fig:minus_1_weighting}).

\begin{lemma}
\label{lemma:minus_1_weight}
Let $\lambda$ be a 2-factor of $H_{a,b,c}$.  Then
\[
\sum_{\mu \in \Psi^{-1}(\lambda)} w_{-1}(\mu) = (-1)^{ab+bc+ca} \cdot  2^{\#\{\text{closed paths in $\lambda$}\}}.
\]
\end{lemma}

\begin{figure}
\caption{The unique configuration of edges which squishes to a doubled edge.}
\label{fig:squish_to_doubled_edge}
\begin{center}
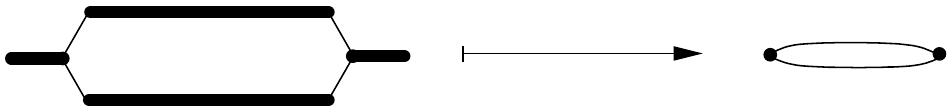
\end{center}
\end{figure}

\pf Let us suppose that $\lambda$ decomposes as the disjoint union of doubled edges $e_1, \ldots, e_n$ and closed loops $\ell_1, \ldots, \ell_m$.  Because the unsquishings $\varphi(e_i)$ and $\varphi(\ell_j)$ are all pairwise non-adjacent in $H_{2a,2b,2c}$,  it is clear that

\begin{equation}
\label{eqn:decomp}
\sum_{\mu \in \Psi^{-1}(\lambda)} w_{-1}(\mu) = 
\prod_{i=1}^n \left( \sum_{\mu \text{ matching on } \varphi(e_i)} \negthickspace \negthickspace \negthickspace w_{-1}(\mu) \right) 
\prod_{j=1}^m \left(  \sum_{\mu \text{ matching on }\varphi(\ell_j)} \negthickspace \negthickspace \negthickspace w_{-1}(\mu) \right).
\end{equation}
Let us first consider $\Psi^{-1}(e_i)$, the preimage of a doubled edge.  There is only one possible matching on $\phi(e_i)$ (see Figure~\ref{fig:squish_to_doubled_edge}); it has weight -1.

\begin{figure}
\caption{Lifting a closed path $\ell_j$ to $\varphi(\ell_j)$, and completing it to a matching.}
\label{fig:lift_walk}
\begin{center}
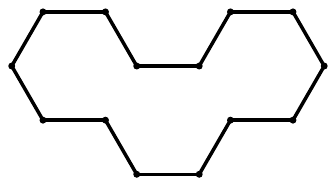
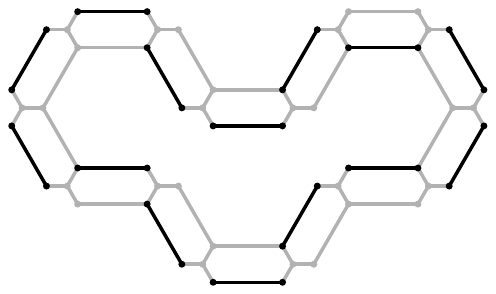
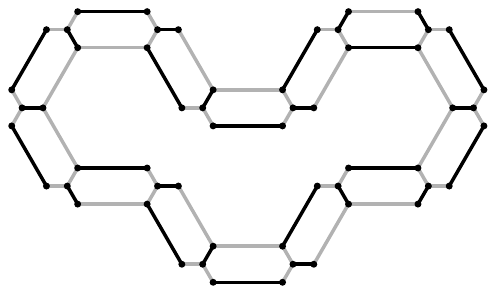
\end{center}
\end{figure}

Next, let us consider $\Psi^{-1}(\ell_j)$, the preimage of a closed loop.  Every term in $\Psi^{-1}(\ell_j)$ is a matching on $\varphi(\ell_j)$ which maps to $\ell_j$ under $\psi$.  We shall ``lift'' $\ell_j$ up to $\varphi(\ell_j)$, so that at each step we choose an edge $e' \in \varphi(\ell_j)$ which maps to $e$ under $\psi$.  If we do this in such a way that no two of the $e'$ are adjacent, then there is a unique way to complete the lifted walk to a matching on $\varphi(\ell_j)$, by choosing exactly one edge of each propeller to be in the matching (see Figure~\ref{fig:lift_walk}).   Note that the $w_{-1}$-weight of the lifting (second image) is the same as the weight of the matching (third image) since all of the propellers have weight 1.

Suppose that we walk along the lifted path in $\varphi(\ell_j)$, in the
counterclockwise (positive) direction, and we step on the edges $e'_1, e'_2,
\cdots, e'_k, e'_1$, in order.  We keep track of the product of the weights of
the edges we have stepped on.

At any point in our walk, there are essentially only four states that we can be in.  We can be standing on an edge weighted 1 or -1, and we can be facing one of two types of propeller: \Yright or \Yleft.  Let us label the states 1 through 4, as follows:

\begin{enumerate}
\item Standing on a weight 1 edge, facing a \Yright propeller
\item Standing on a weight -1 edge, facing a \Yright propeller
\item Standing on a weight 1 edge, facing a \Yleft propeller
\item Standing on a weight -1 edge, facing a \Yleft propeller
\end{enumerate}
 
To get from one edge to another in the path, we must turn either left or right at the propeller we are facing, so we shall think of our walk as a sequence of Rs and Ls, read \emph{right-to-left}.   For example, the walk in Figure~\ref{fig:lift_walk}, starting at the upper leftmost horizontal edge and proceeding counterclockwise, is $LRRLLLLRLLRLLL$.  

Our trick is to interpret $R$ and $L$ as the following \emph{state--transition matrices}:
\begin{align*}
L &= \left[
\begin{matrix}
0	& 0	& 1	& 1 	\\
0	& 0	& 0	& -1	\\
1	& 0	& 0	& 0 	\\
-1	& -1	& 0	& 0 	
\end{matrix} \right]
& R &= \left[
\begin{matrix} 
0	& 0	& 1	& 0 	\\
0	& 0	& -1	& -1 	\\
1	& 1	& 0	& 0 	\\
0	& -1	& 0	& 0 	
\end{matrix} \right]
\end{align*}

These matrices encode which states can be reached from which, and keep track of the total weight of all paths leading from one state to another.  For example, at the upper-left most horizontal edge in Figure~\ref{fig:lift_walk}, we are in state 4, and we are about to turn to the left.  We could go into state 1, picking up a weight of 1 (as in the diagram); alternately, we could go into state 2, picking up a weight of -1.  This is modeled by our state-transition matrices, since
\[
L \cdot \left[\begin{matrix} 1 \\ 0 \\ 0 \\ 0  \end{matrix}\right]
=
\left[
\begin{matrix}
0	& 0	& 1	& 1 	\\
0	& 0	& 0	& -1	\\
1	& 0	& 0	& 0 	\\
-1	& -1	& 0	& 0 	
\end{matrix} \right] \cdot
\left[\begin{matrix} 1 \\ 0 \\ 0 \\ 0  \end{matrix}\right]
= \left[\begin{matrix} 0 \\ 0 \\ 1 \\ -1  \end{matrix}\right]
\]
By the same token, the coefficients of the matrix $LR$ are the sums of the weights of all paths which start in state $i$ , turn first right and then left, and end in state $j$.  If we want to know the sum of the weights of the lifts of the closed loop in Figure~\ref{fig:lift_walk}, we should compute the matrix product $LRRLLLLRLLRLLL$.  We must end in the same state as we start since our path is closed; furthermore, our starting state must be either 3 or 4, so the total number of paths is the sum of the (3,3) and (4,4) entries of this matrix product.

Now, here is the punchline:  It is easy to check that $L=R^{-1}$.   In any closed loop on the hexagonal lattice, there must be 6 more left turns than right turns, so the combined weights of all lifts of any closed loop is the sum of the $(3,3)$ and $(4,4)$ entries of $L^6$.  Now, it is also easy to check that $L^6 = -I$, so this combined weight is -2, regardless of the shape or orientation of the loop $\ell_j$.

Going back to Equation~\ref{eqn:decomp}, and recalling that there are $i$ doubled edges and $j$ closed loops in the 2-factor $\lambda$, we can now state that
\[
\sum_{\mu \in \Psi^{1}(\lambda)} w_{-1}(\mu) = (-1)^i \cdot (-2^j) = (-1)^{i+j} \cdot 2^j
\]
Observe that $i+j$ is the total number of connected components in $\lambda$, so the proof is complete modulo the following lemma. $\hfill \square$  

\begin{lemma}
\label{lemma:component_parity}
Let $\lambda$ be a 2-factor of $H_{a,b,c}$.  Let $C(\lambda)$ denote the number of connected components of $\lambda$.  Then $C(\lambda) \equiv ab+bc+ca \pmod{2}$.
\end{lemma}

\pf Decompose $\lambda$ into two matchings, $\lambda_1$ and $\lambda_2$.  Let us consider $\lambda_!$ for the moment.  

\begin{figure}
\caption{A $\tau$-transformaion}
\label{fig:local_matching_change}
\begin{center}
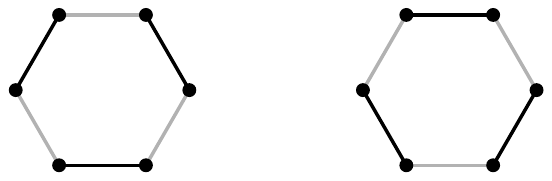
\end{center}
\end{figure}

When one removes (or adds) a cube to the 3D diagram corresponding to $\lambda_1$, one obtains a new matching.  Graph-theoretically, this is equivalent to performing the change depicted in Figure~\ref{fig:local_matching_change}.  Let us call such an operation a \emph{$\tau$-transformation}.  Since it is possible to change any 3D diagram into any other by adding and removing boxes one at a time, repeated $\tau$-transformations will allow us to change any matching into any other matching.

\begin{figure}
\caption{Changing one 2-factor into another with $\tau$-transformations}
\label{fig:local_2_factor_change}
\begin{center}
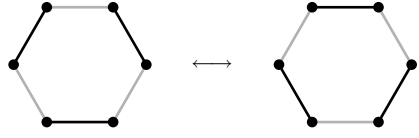
\end{center}
\end{figure}

Let us alter the matchings $\lambda_1$ and $\lambda_2$ in this fashion, so that after finitely many steps, we have obtained two copies of the empty 3D Young diagram.   At each stage, superimposing the two matchings gives a 2-factor, so we may think of $\tau$-transformations as acting on the set of 2-factors as well (see Figure~\ref{fig:local_2_factor_change}).

\begin{figure}
\caption{The two possible connectivities of $\lambda_1^{\circ} \cup \lambda_2$, and their behaviour under $\tau$-transformations}
\label{fig:delete_a_hexagon}
\begin{center}
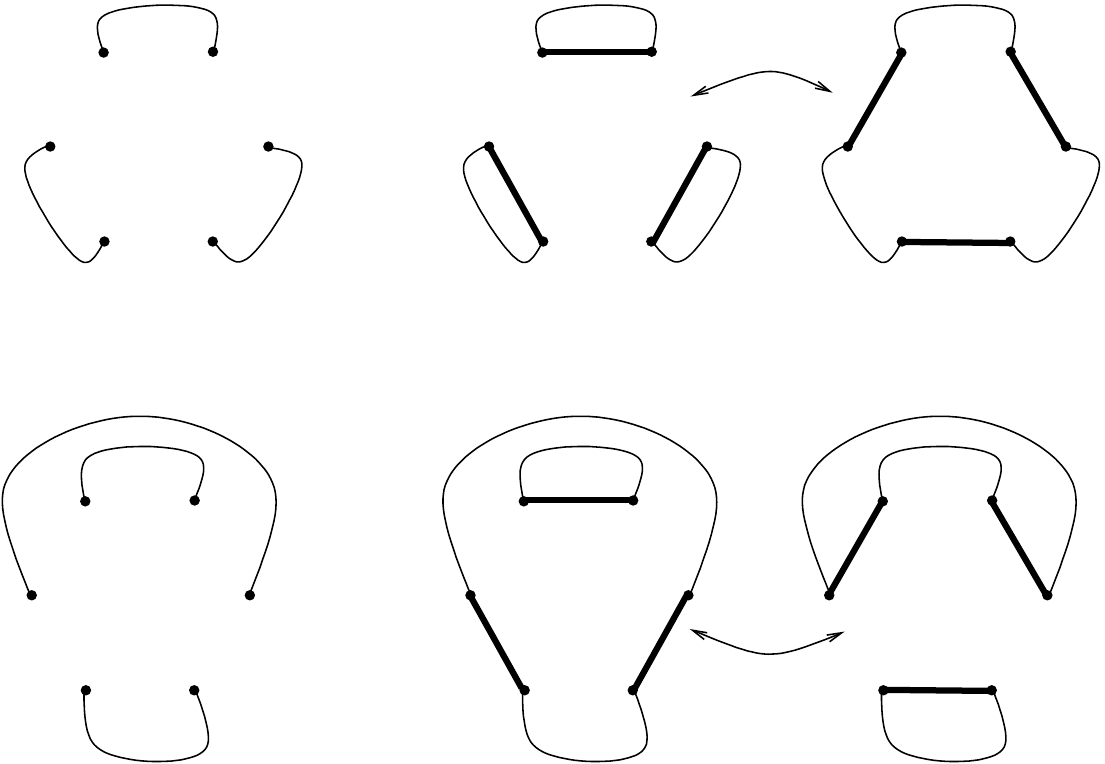
\end{center}
\end{figure}

Next, we show that $\tau$-transformations leave $C(\lambda) $ invariant mod 2.
Suppose we are applying the $\tau$ transformation to $\lambda_1$, and we are affecting the edges around a particular hexagon $H$.  Remove all edges of $H$ from $\lambda_1$ to obtain a new subgraph $\lambda_1^{\circ}$ of $H_{a,b,c}$.  Now superimpose $\lambda_1^{\circ}$ and $\lambda_2$.  The resulting graph has degree 2 everywhere except at vertices of $H$, all of which have degree 1.  As a result, every vertex of $H$ is connected to precisely one other vertex of $H$ by some path (possibly of length one).   Up to symmetry, there are two possible ways in which this may occur (see Figure~\ref{fig:delete_a_hexagon}).  When we add the edges of $\lambda_1 \cap H$ back into the graph and perform the $\tau$-transformation, it is easy to see that the parity of the number of connected components is preserved.  In the upper case, there are an odd number of path components; in the lower case, there are an even number. 

Up to this point, we have shown that $C(\lambda)$ is an invariant of $H_{a,b,c}$.  To calculate it, we compute the number of edges in the matching corresponding to the empty 3D Young diagram.  This is easy: there are $ab$ horizontal edges, $ac$ edges of slope $\sqrt{3}$, and $bc$ edges of slope $-\sqrt{3}$, and so there are $ab+bc+ca$ edges in all.
$\hfill \square$

\begin{lemma}
\label{lemma:p111_lemma}
If $\lambda$ is a 2-factor on $H_{a,b,c}$ and $\mu \in \Psi^{-1}(\lambda)$, then
$w_{p,1,1,1}(\lambda) = w_p(\mu)$.
\end{lemma}

\begin{figure}
\caption{The weightings $w_{p,1,1,1}$ on $H_{2a,2b,2c}$, compared with the weighting $w_p$ on $H_{a,b,c}$, where $t=p^3$.}
\label{fig:p111_weighting}
\begin{center}
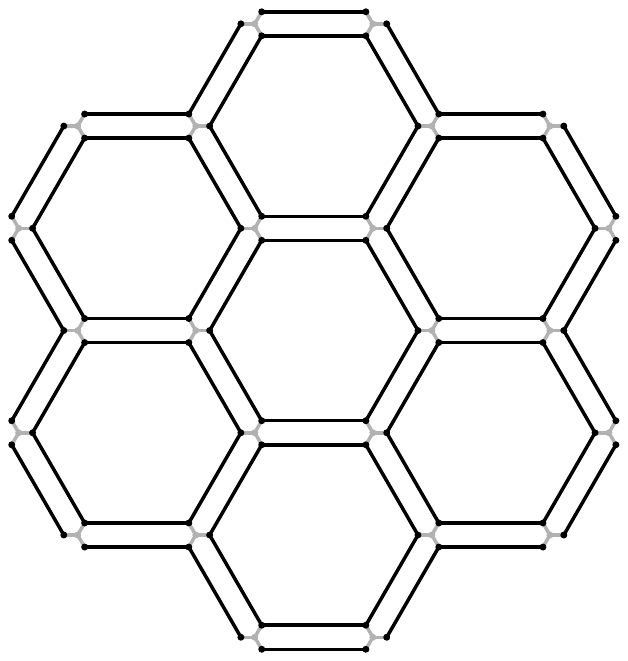
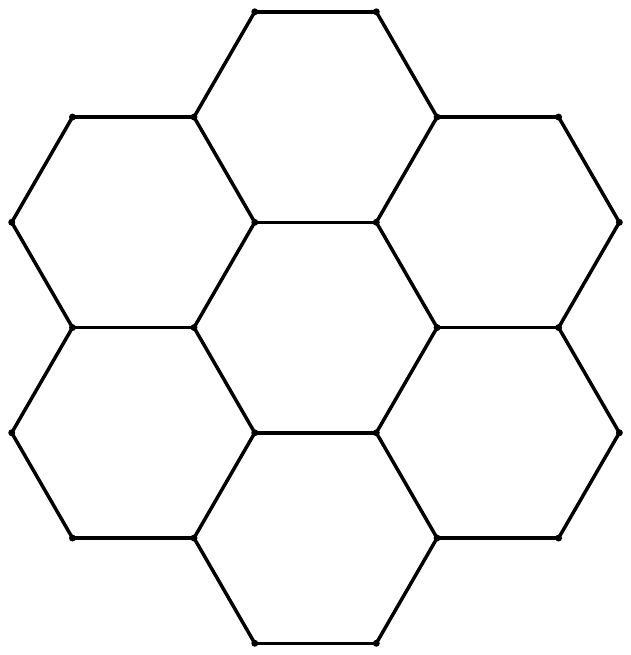
\end{center}
\end{figure}

\pf Comparing the weightings $w_{p,1,1,1}$ and $w_p$ in Figure~\ref{fig:p111_weighting}, it is clear that if $e$ is an edge of $H_{a,b,c}$ and $e'$ is either of the two edges in $\psi^{-1}(e)$, then $w_{p,1,1,1}(e') = w_{p}(e)$.  The statement of the lemma follows immediately. $\hfill \square$.

\section{A specialization of the $\bbZ_2 \times \bbZ_2$ partition function}

As before, $w_{p, q, r, s}$ denotes the $\bbZ_2 \times \bbZ_2$ weighting, and
$w_p =
w_{p,p,p,p}$ denotes the monochromatic weighting.  We will be working with
various specializations of these weightings, and it will be convenient to
write, for example, $w_{-p, -1, -1, -1}$ for the specialization $p \mapsto -p$,
$q,r,s \mapsto -1$.   /usr/X11R6/bin/gvim

Although our tool of choice is weighted matchings on hexagon meshes, our object of study is the $\bbZ_2 \times \bbZ_2$ weighting for 3D Young diagrams.  We shall use $w_{p,q,r,s}(\pi)$ to denote the $\bbZ_2 \times \bbZ_2$--weight of the 3D diagram $\pi$, and $w_{p}$ to denote its monochromatic weight.  Recall that if $\lambda$ is the matching corresponding to $\pi$, then 
\[
w_{p,q,r,s}(\pi) = \frac{w_{p,q,r,s}(\lambda)}{w_{p,q,r,s}(\text{empty 3D diagram})}.
\]
and we can normalize the monochromatic weight in the same fashion.
\begin{definition} The $\bbZ_2 \times \bbZ_2$  partition function for $H_{a,b,c}$ is
\[
Z^{a,b,c}(p,q,r,s) =
\sum_{\substack{ \pi \text{ 3d diagram} \\ \text{inside }a \times b \times c \text{ box} }}  \!\!\!\!\!\!\!\!w_{p, q, r, s}(\pi).
\]
The monochromatic partition function for $H_{a,b,c}$ is
\[
Z^{a,b,c}(p) = 
\sum_{\substack{ \pi \text{ 3d diagram} \\ \text{inside }a \times b \times c \text{ box} }} \!\!\!\!\!\!\!\!w_{p}(\pi)
\]
\end{definition}

Note that $\lim_{a,b,c \rightarrow \infty}Z^{a,b,c} = \lim_{a,b,c \rightarrow \infty}Z^{2a,2b,2c} = Z_{\bbZ_2 \times \bbZ_2}$.  Thus the following implies \eqref{eqn:-1_specialization}:

\begin{theorem}
\[Z^{2a,2b,2c}(p,-1,-1,-1) =  \left(Z^{a,b,c}(-p)\right)^2. \]
\end{theorem}

\pf To begin with, let us work with matchings.  Observe that for any matching $\mu$, we have $w_{p,-1,-1,-1}(\mu) = w_{-p,1,1,1}(\mu) \cdot w_{-1,-1,-1,-1}(\mu)$.   Therefore,

\begin{align}
\nonumber
\sum_{\substack{\mu \text{ matching} \\ \text{on } H_{2a,2b,2c}}}\!\!\!\! w_{p,-1,-1,-1}(\mu) 
&= \sum_{\substack{\lambda \text{ 2-factor} \\\text{on }H_{a,b,c}}} \sum_{\mu \in \Psi^{-1}(\lambda)} w_{-p, 1, 1, 1}(\mu) w_{-1, -1, -1, -1}(\mu) \\
\label{eqn:first_guy}
&= \sum_{\substack{\lambda \text{ 2-factor} \\\text{on }H_{a,b,c}}} w_{-p}(\lambda) \sum_{\mu \in \Psi^{-1}(\lambda)} w_{-1, -1, -1, -1}(\mu) \\
\label{eqn:second_guy}
&= (-1)^{ab+bc+ca}\sum_{\substack{\lambda \text{ 2-factor} \\\text{on }H_{a,b,c}}} 2^{\#\{\text{closed loops in }\lambda\}}w_{-p}(\lambda) \\
\label{eqn:third_guy}
&= (-1)^{ab+bc+ca} \sum_{\substack{\lambda, \eta \text{ 1-factors}\\ on H_{a,b,c} }}w_{-p}(\lambda)w_{-p}(\eta)\\
\nonumber
&= (-1)^{ab+bc+ca}\left(
\sum_{\substack{\mu \text{ matching} \\ \text{on } H_{a,b,c}}}\!\!\!\! w_{-p}(\mu)
\right)^2
\end{align}
In the above sequence of steps, Equation~(\ref{eqn:first_guy}) uses Lemma~\ref{lemma:p111_lemma}, Equation~(\ref{eqn:second_guy}) uses Lemma~\ref{lemma:minus_1_weight}.

Now, let us normalize both sides by dividing by the $w_{p,-1,-1,-1}$ weight of the matching associated to the $(2a,2b,2c)$ empty 3D diagram.  Let us call this matching $M$.  Call the $(a,b,c)$ empty 3D diagram $N$.  

On the left side, we get $Z^{a,b,c}(p,q,r,s)$.  For the right side, we observe that 
\begin{align*}
w_{p,-1,-1,-1}(M) &= w_{-1,-1,-1,-1}(M) \cdot w_{-p,1,1,1}(M)\\
&= (-1)^{ab+bc+ca} \cdot (w_{-p}(N))^2,
\end{align*}
again using Lemma~\ref{lemma:p111_lemma}.  Therefore, the right-hand side normalizes to 
\[
\left(Z^{a,b,c}(-p)\right)^2.
\]
and we are done. $\hfill \square$

\section{Conclusions}

It is obvious to ask whether this method might be used for its intended purpose, namely to prove Equation~\ref{eqn:z2z2_function}.  Unfortunately, the answer is no.  Our method is clearly agnostic as to the large-scale shape of the graph; it works on any suitable subset of the hexagon lattice.  However, computer calculations show that the analogue of~\eqref{eqn:z2z2_function} on $H_{2a,2b,2c}$ does not admit a product formula; rather, if one counts 3D Young diagrams fitting inside a $2a \times 2b \times 2c$ box, according to the $\bbZ_2 \times \bbZ_2$ weighting, one typically gets $(1-p) \times (\text{a large irreducible polynomial})$; the product formula only appears as $a,b,c \rightarrow \infty$.

It should be possible to apply this method to other graphs and other weights.  We have not investigated this possibility in any great depth, but in order for the method to work without serious modification, it seems that the graph ought to have propeller-like vertices, at which several hexagonal faces meet; furthermore, one might expect the graph to have some degree of self-similarity at different scales.

\section{Acknowledgements}

We would like to thank Jim Bryan for editing and proofreading, and Richard Kenyon for the idea of the matrices $L$ and $R$ in section 5. 

\bibliographystyle{amsplain}
\bibliography{squish}

\end{document}